\documentstyle[12pt,amssymb,psfig,epsfig,color]{amsart}
\def\from#1\to{\mathpunct:#1\rightarrow}

\newtheorem{thm}{Theorem}[section]
\newtheorem{cor}[thm]{Corollary}
\newtheorem{lem}[thm]{Lemma}

\theoremstyle{remark}

\theoremstyle{definition}

\numberwithin{equation}{section}
\numberwithin{figure}{section}

\def\note#1{}

\newcommand{\di}{\partial}

\newcommand{\ra}{\rightarrow}

\def\ssk{\smallskip}
\def\msk{\medskip}

\def\nin{\noindent}

\def\sm{\smallsetminus}

\renewcommand{\mod}{\operatorname{mod}}

\newcommand{\tl}{\tilde}

\newcommand{\Time}{\operatorname{Time}}

\newcommand{\orb}{\operatorname{orb}}

\newcommand{\eps}{{\varepsilon}}

\newcommand{\De}{{\Delta}}
\newcommand{\de}{{\delta}}

\newcommand{\La}{{\Lambda}}

\newcommand{\Om}{{\Omega}}

\newcommand{\OO}{{\cal O}}

\newcommand{\RR}{{\cal R}}

\newcommand{\TT}{{\cal T}}

\newcommand{\C}{{\Bbb C}}

\newcommand{\N}{{\Bbb N}}

\newcommand{\W}{{\Bbb W}}

\def\B0{{\mathbf{0}}}

\catcode`\@=12

\def\Empty{}
\newcommand\oplabel[1]{
  \def\OpArg{#1} \ifx \OpArg\Empty {} \else
  	\label{#1}
  \fi}
		
\newcommand{\comm}[1]{}
\newcommand{\comment}[1]{}

\begin{document}

\bigskip\bigskip

\title[Higher degree polynomials]{Local connectivity of Julia sets\\ for unicritical polynomials}
\author {Jeremy Kahn and Mikhail Lyubich}
\date{April 29, 2005}

\thispagestyle{empty} 
\def\IMSmarkvadjust{0 pt}
\def\IMSmarkhadjust{0 pt}
\def\IMSmarkhpadding{0 pt}
\def\IMSpubltext{Published in modified form:}
\def\SBIMSMark#1#2#3{
 \font\SBF=cmss10 at 10 true pt
 \font\SBI=cmssi10 at 10 true pt
 \setbox0=\hbox{\SBF \hbox to \IMSmarkhpadding{\relax}
                Stony Brook IMS Preprint \##1}
 \setbox2=\hbox to \wd0{\hfil \SBI #2}
 \setbox4=\hbox to \wd0{\hfil \SBI #3}
 \setbox6=\hbox to \wd0{\hss
             \vbox{\hsize=\wd0 \parskip=0pt \baselineskip=10 true pt
                   \copy0 \break%
                   \copy2 \break%
                   \copy4 \break}}
 \dimen0=\ht6   \advance\dimen0 by \vsize \advance\dimen0 by 8 true pt
                \advance\dimen0 by -\pagetotal
	        \advance\dimen0 by \IMSmarkvadjust
 \dimen2=\hsize \advance\dimen2 by .25 true in
	        \advance\dimen2 by \IMSmarkhadjust

%
%
  \openin2=publishd.tex
  \ifeof2\setbox0=\hbox to 0pt{}
  \else 
     \setbox0=\hbox to 3.1 true in{
                \vbox to \ht6{\hsize=3 true in \parskip=0pt  \noindent  
                {\SBI \IMSpubltext}\hfil\break
                \input publishd.tex 
                \vfill}}
  \fi
  \closein2
  \ht0=0pt \dp0=0pt
 \ht6=0pt \dp6=0pt
 \setbox8=\vbox to \dimen0{\vfill \hbox to \dimen2{\copy0 \hss \copy6}}
 \ht8=0pt \dp8=0pt \wd8=0pt
 \copy8
 \message{*** Stony Brook IMS Preprint #1, #2. #3 ***}
}

\SBIMSMark{2005/03}{April 2005}{}

\begin{abstract} 
  We prove that the Julia set $J(f)$ of 
at most finitely renormalizable unicritical polynomial  $f:z\mapsto  z^d+c$  with all periodic points repelling 
is locally connected. (For $d=2$ it was proved by Yoccoz around 1990.)
It follows from {\it a priori} bounds in a modified Principle Nest of puzzle pieces. 
The proof of {\it a priori} bounds  makes use of new analytic tools developed in \cite{covering lemma}
that give control of moduli of annuli under maps of high degree. 
\end{abstract}

\setcounter{tocdepth}{1}
 
\maketitle


\section{Introduction}\label{intro}

\subsection{Statement of  the results}  

About 15 years ago Yoccoz proved that the Julia set of at most finitely many renormalizable quadratic polynomial
$f: z\mapsto z^2+c$ with all periodic points repelling is locally connected (see \cite{H,M}).
In this paper, we generalize this result to higher degree unicritical polynomials: 

\proclaim Theorem A. 
The Julia set $J(f)$ of 
at most finitely renormalizable unicritical polynomial  $f: z\mapsto  z^d+c$  with all periodic points repelling 
is locally connected.

This result follows from {\it a priori} bounds in an appropriate 
``Modified Principle Nest''  of puzzle pieces,
$$
  E^0\Supset E^1\Supset\dots \ni 0: 
$$

\proclaim Theorem B. The (modified) principal moduli stay away from zero:
$$
    \mod(E^{i-1}\sm E^i)\geq \mu > 0.
$$

These {\it a priori} bounds imply that the puzzle pieces $E^i$ shrink to the critical point,
which yields Theorem A by a standard argument.%

\subsection{Techniques}

As usual in holomorphic dynamics, our proof has two sides: combinatorial and analytic.
Our combinatorial tool is a refined  Principal Nest techniques of \cite{puzzle}, 
while the analytic tool is a recently established
Quasi-Invariance Law (Covering Lemma) in conformal geometry \cite{covering lemma}.  
Let us briefly comment on both sides. 

The puzzle machinery was introduced to holomorphic dynamics by Branner and Hubbard \cite{BH}
(in the context of cubic polynomials with one escaping critical point) and
Yoccoz \cite{H,M} (in the context of quadratic polynomials). 
The idea is to tile shrinking neighborhoods
of the Julia set into topological disks called {\it puzzle pieces}, 
and to translate the dynamics on $J(f)$ to the combinatorics of these tilings.

An efficient way to describe these combinatorics is given by
the {\it Principal nest} of puzzle pieces around the origin,
$V^0\supset V^1\supset\dots V^n\dots\ni 0$, which is inductively constructed so that 
the first return maps $f^{n_i}: V^{i}\ra V^{i-1}$ are unicritical branched coverings
 \cite{puzzle}. 
It turns out that this nest is not quite suitable for our purposes, 
so we modify it slightly to obtain  
a {\it dynasty of kingdom map}, see \S \ref{principal nest}. 

We then observe that since the return times in the dynasty grow exponentially,
one can send some puzzle piece $E^{i-1}$ 
to the top level by an appropriate composition $\Psi$ of the kingdom maps, 
while the next puzzle piece, $E^i$, will go at most five levels up ({\it time inequality}).
Thus, the map $\Psi|\, E^i$ has a bounded degree, which 
puts us in a position to apply the analytic techniques of \cite{covering lemma}.
 
The puzzle bears  complete information about the Julia set only if
the puzzle pieces shrink to points, so this is a key geometric issue of the theory.
To handle this issue,  Branner \& Hubbard and Yoccoz made use of the 
Series Law from conformal geometry.%
\footnote{also called  the   Gr\"otzsch Inequality.}
It was immediately realized, however, that this method would not work for higher degree polynomials,
so that in the higher degree case the problem has remained open since then.

A new analytic tool that we exploit is a Covering Lemma (Quasi-Invariance Law) in conformal geometry \cite{covering lemma}
which roughly asserts that given a branched covering $g: U\ra V$ of degree $N$ which restricts to a branched covering
$g: A\ra B$ of degree $d$ such that $\mod (U\sm A) $ is small (depending on $N$), 
then, under a certain ``Collar Assumption'',  $ \mod(V\sm B)$  is  comparable with  $d^2 \mod(U\sm A)$  
(independently of $N$) -- see \S \ref{covering lemma sec} for the precise statement.    

The Covering Lemma allows us to transfer moduli information from deep levels of the dynasty to 
shallow ones, and to argue that if on some deep levels the moduli are small, then they must be even 
smaller on  shallow ones. This certainly implies that, in fact, the moduli can never  be too small (Theorem B).

\ssk
Note in conclusion that
for real $c$, Theorem A  was proved before by Levin \& van Strien \cite{LS}.
The method  used in \cite{LS} exploited real symmetry in a substantial way. 

\ssk In the forthcoming notes (joint with A. Avila and W. Shen)
      our {\it a priori} bounds will be used to prove {\it rigidity} of the 
      unicritical polynomials under consideration.

\subsection{Terminology and Notation}   

  A {\it topological disk} means a simply connected  domain in $\C$.
 
  We let $\orb(z)\equiv \orb_g(z)= (g^n z)_{n=0}^\infty$ be the {\it orbit} of $z$ under a map $g$.

  Given a map $g: U\ra V$ and a domain $D\subset V$, 
components of $g^{-1}(D)$ are called {\it pullbacks}  of $D$ under $g$. 
Given a connected set $X\subset g^{-1}(D)$, we let $g^{-1} (D)|X$ be the pullback of $D$
containing $X$.    

Given a subset $W\subset V$, the {\it first landing map} $H$ to $W$ is defined (on the set of points $z$
whose orbits intersect $W$) as follows:  $H(z)= f^l z$, where $l\geq 0$ is the first moment for which
$f^l z\in W$.

We say that a map $g: U\ra V$ is {\it unicritical} if it has one critical point
(of arbitrary local degree)

\msk {\bf Acknowledgment.}
 We thank Artur Avila for careful reading the manuscript and making a number of useful comments.
We also thank all the Foundations that have supported this work: 
the Guggenheim Fellowship, Clay Mathematics Institute, NSF, and NSERC.

\section{Modified Principal Nest}\label{principal nest}

\subsection{Generalized polynomial-like maps}\label{GPL}

 A {\it generalized polynomial-like map} (GPL map) is a holomorphic map
$g\from \cup W_i \to V$, where $V \subset \C$ is a topological disk and 
$W_i\Subset V$ are topological disks with disjoint closures such that the restrictions
$g:W_i\ra V$ are branched coverings, and moreover, all but finitely many of them have degree one.

\ssk {\it Remark.} To prove Theorem B in full generality, we need to allow infinitely many disks
$W_i$. However, in the ``persistently recurrent'' case that interest us most
it is enough to consider GPL maps defined on finitely many disks $W_i$.
\ssk

\comm{********
 A {\it generalized polynomial-like map} (GPL map) is a proper holomorphic map
$g\from \W \to V$, where $V \subset \C$ is a topological disk and $\W \Subset V$ is full.
Then $\W = \bigcup_{i=0}^k W_i$, where the domains $W_i$ are topological disks
 and the map $g: W_i \ra V$ is a branched covering. 
***************}

We let $K_g  = \bigcap_{n=0}^\infty g^{-n} V$ be the set of points of $V$ on which $g$ is infinitely iterable
(the ``filled Julia set'').

%

A GPL map $g$ is called {\it unicritical} if it has a single critical point.
{\it In what follows we will consider only unicritical  GPL maps, and we will always put its critical point at $0$.} 
Let $d$ be  the local degree of $g$ near 0. We let   $W_0\equiv W$ be the ``central domain'', that is, the one
containing~0.  

The {\it postcritical set} $\OO_g$ of a (unicritical) GPL map is the closure of the orbit $\{g^n 0\}_{n=0}^\infty$.

\msk
  {\it Puzzle pieces} of depth $n$ of a GPL map $g$ are components of $g^{-n}(V)$. Puzzle pieces containing
$0$ are called {\it critical}. 

If the critical point returns to some critical puzzle piece $A$,
then  the first return map $h$ to $A$ is also GPL. Let   $\cup B_i$ be its domain of definition.  
Restricting $h$ to the union of those components $B_i$ that intersect the postcritical set,
we obtain a GPL map called the {\it generalized renormalization} 
$r_A(g)$ of $g$ on $A$. 

If we do not specify the domain $A$ of the generalized renormalization,
then it is assumed to be $W$, so $r(g)\equiv r_W(g)$.


\subsection{Dynasty of kingdoms}\label{dynasty}
Let us introduce a modified notion of (unicritical) GPL map called a kingdom map. 

Let us consider three topological disks,  $W\supset U\Supset  A\ni 0$,
called the {\it kingdom domain}, the {\it castle},  and the {\it king} respectively. 
Let $D_j\Subset W\sm\bar A$ be a family of topological disks (`` king's subjects'')
such that $\bar D_j\cap \di U=\emptyset$.
Finally, let  $M_k\Subset U\sm\bar A$ be another family of topological disks 
(``king's men'').
A map
$$
     G: A \cup D_j\cup M_k \ra W
$$
is called a {\it kingdom map} (of {\it local degree} $d$) if
\begin{itemize}
 \item  The closures $\bar A$, $\bar D_j$ and $\bar M_k$ are pairwise disjoint; 
 \item $G: A\ra W$ is  a $d$-to-$1$ branched covering ramified only at 0;
 \item Each $G: D_j\ra W$  is a biholomorphic isomorphism; 
 \item Each $G: M_k\ra U$ is is a biholomorphic isomorphism. 
\end{itemize}
We let $\OO_G$ be the postcritical set of the kingdom map $G$.

\begin{figure}[htbp]
\begin{center}
\input{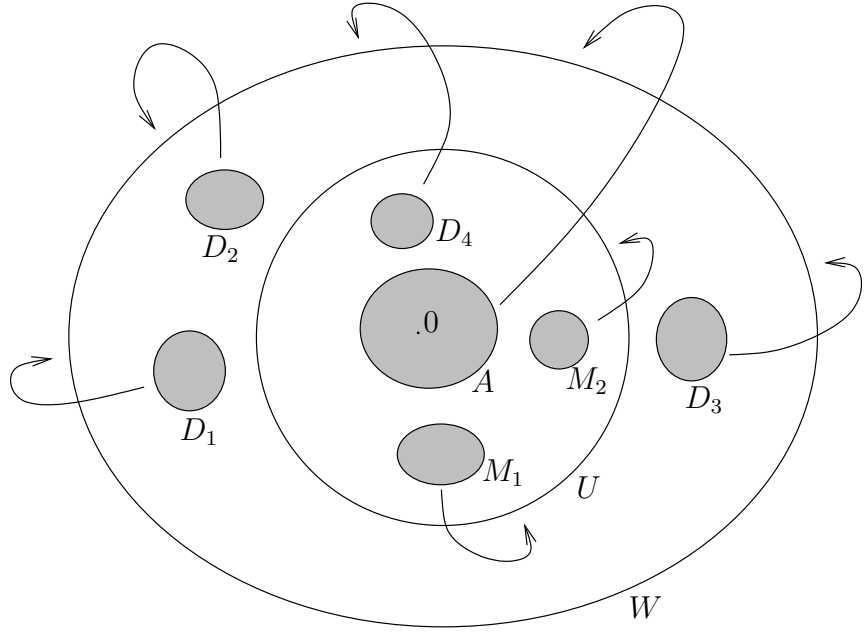}
\caption{Kingdom map}
\end{center}
\end{figure}


When $U=W$, kingdom maps become GPL maps.

\msk
Let us now consider a (unicritical)  GPL map $g: \cup W_i\ra V$, $W\equiv W_0$. 
Let us define  the {\it kingdom renormalization} $G=R(g)$ of $g$
whose result will be a kingdom map $G$.

 If $g(0)\in W$ then we say that the {\it central return} occurs.
If $g^k (0)\in W$ for $k=0, \dots, N-1$ but $g^N(0)\not\in W$, 
then we have a nest of topological disks
\begin{equation}\label{cascade}
    V\equiv \Om_0\Supset W\equiv  \Om_1\Supset\dots \Supset \Om_N\equiv U
\end{equation}
such that $g: \Om_{k+1}\ra \Om_k$ is a unicritical branched covering of degree $d$ 
and $g(0)\in \Om_{N-1}\sm \Om_N$.   
This nest is called a {\it central cascade (of length $N$)}. 
Note that the non-central-return event corresponds to the cascade of length 1. 

In  the kingdom renormalization $Rg$,  $W$ will be the kingdom domain and $U$ will be the  castle.

Let us consider the first return map to $W=\Om_1$:
$$
   h : X_0 \cup \bigcup_{i>0} X_i  \ra W, 
$$
where $X_0\ni 0$ (in case when $N\geq 2$, $X_0= \Om_2$). 

Let us consider the domain $X_s$, $s>0$,  containing  $g^{N-1}(0)$. 
Then  the pullback $A  = g^{-(N-1)}(X_s)|\, 0$ is the  king of $Rg$.  
The kingdom map $G$ on $A$ is defined as $ h \circ g^{N-1}: A \ra W$.
Notice that it is a unicritical $d$-to-$1$ branched covering. 

Let us define {\it king's subjects} $D_j$ as non-critical pullbacks of the domains $X_i$ ($i\not= 0$) under the maps 
$$
   g^{k-1}: \Om_k\ra W,\quad k=1,2,\dots, N,
$$
that intersect the postcritical set. 
Thus,  each subject $D_j$ is  univalently mapped  onto some $X_i$, $i>0$,
by an appropriate map $g^{k-1}: \Om_k \ra W$, $k\in [1,N]$.
On this subject let us define the kingdom  map $G: D_j\ra W$ as $G|D_j = h\circ g^{k-1}|D_j $.
Obviously, it is a biholomorphic isomorphism. 


Finally, let us define  king's men $M_k$ 
as the pullbacks of $U$ under  $g: U\ra g(U)\supset U$ that intersect the postcritical set.
There are at most $d$ king's men, and $g$ univalently maps each of them onto $U$. 
Let $G|M_k= g|M_k$. 

Thus, we have defined the desired kingdom renormalization 
$$
   G=R(g) : A \cup D_j\cup M_k \ra W.
$$ 

Let $g^N(0)\in W_j$, $j>0$. For $G=R(g)$,  
let us define the {\it king's apartment} $\De$ as $g^{-N}(W_j)|\, 0$. 
Then $A\Subset \De\subset U$
and the map $g^{N+1} : \De\ra V$ is a unicritical $d$-to-1 branched covering. 
This creates a collar $\De\sm A$ around the king.  

\ssk
{\it Remark.} If $N=1$ (i.e., the non-central return occurs under $g$), then the kingdom renormalization
$G=R(g)$ coincides with the generalized renormalization defined in \S \ref{GPL}. 


\msk
 Given a kingdom map $G$, let us define its renormalization $g=r(G)$ as the first return map
$g: \cup B_i\ra A$  to the king $A$ restricted to those domains $B_i$ that  intersect the postcritical set $\OO_G$. 
It is a unicritical GPL map.


Beginning with some GPL map $g\equiv g_0$, 
we construct in the above way  a {\it dynasty} of kingdoms, that is, 
 a sequence  $(g_n, G_n)$ such that $g_n$ is a GPL map, $G_n$ is a kingdom map,
$G_n= R(g_n)$ and $g_{n+1} = r(G_n)$. This dynasty terminates if and only if:

\ssk\nin $\bullet$ The map $g$ is combinatorially non-recurrent, that is, 
   the critical point does not return to some critical puzzle piece; or 

\ssk\nin $\bullet$
 Some map $g_n$ has an infinite central cascade, i.e.,
 it is a Douady-Hubbard polynomial-like map \cite{DH}
with non-escaping critical point.
In this case $g$ is called {\it renormalizable in the sense of Douady and Hubbard}.  

\ssk
When we consider a dynasty of kingdoms $(g_n, G_n)$,  the associated domains will be marked with
superscript $n$ (e.g., $V^n$, $W^n$, etc.) 
However, we usually skip the label when we are concerned with a single kingdom.

\ssk {\it Remark.} 
  It is easy to see that the maps $g_n$ coincide with the generalized renormalizations of $g$ on
domains $V^n$ as defined in \S \ref{GPL}, i.e., $g_n= r_{V^n} (g)$.


The nest 
$$
  V^0\supset W^0\supset \dots\supset    W^{n-1} \supset   V^n \supset W^n \supset V^{n+1}\supset\dots
$$
is called a {\it Modified Principal Nest}. 
Sometimes it is convenient to relabel it  in a uniform way: 
\begin{equation}\label{E}
  E^0\supset E^1\supset  \dots \supset  E^i \supset E^{i+1} \supset\dots,
\end{equation}
so that $V^n=E^{2n}$, $W^n= E^{2n+1}$. 
The consecutive $E$-domains are dynamically related: $E^{i-1}= \psi_i(E^i)$, 
where $\psi_i$ is a unicritical $d$-to-1 branched covering which is an appropriate iterate of $g$.

\subsection{First king}\label{sec: first king}
We will describe in this section how to  associate to a unicritical polynomial $f:z\mapsto z^d+c $ 
(or, more generally, polynomial-like map) a  dynasty of kingdom maps.   
Our standing assumption is that {\it the Julia set $J_f$ is connected and all periodic points of $f$ are
repelling}.  Then $f$ has $d-1$ {\it non-dividing} fixed points $\beta_i$ (landing points of the external rays%
\footnote{In the case of polynomial-like map, external rays are defined by means of straightening.}
with angles $2\pi/(d-1)$), and one {\it dividing} fixed point $\alpha$. There are $q>1$ external rays $\RR_i^0$ landing at 
$\alpha$ which are cyclically permuted by the dynamics, see \cite{M-rays}.

Let us select some equipotential $E^0$; it bounds some topological disk $Q^0$.  
The rays $\RR_i^0$ divide $Q^0$ into $q$ disks $Y^0_i$ called the {\it Yoccoz puzzle pieces} of depth 0.  
Let $Y^0\equiv Y^0_0$ stand for the critical puzzle piece, i.e., the one containing 0. 

The equipotential  $E^1=f^{-1}E^0$ bounds some topological disk $Q^1$.
Let us consider  $dq$  rays of  $f^{-1} (\cup\RR_i^0)$. 
They divide $Q^1$ into $(q-1)d+1$ topological disks called Yoccoz puzzle pieces of depth 1. 
Let $Y^1$ stand for the critical puzzle piece of depth 1. 
There are also $q-1$  puzzle pieces $Y^1_i$ of depth 1 contained in the corresponding off-critical pieces of depth 0.
All other puzzle pieces of depth 1 will be denoted $Z_j^1$. They are attached to the $f$-preimages of $\alpha$
that are different from $\alpha$ itself. 

  The map $f$ is called   {\it satellite renormalizable}
(or, {\it immediately renormalizable})  if
 $$
           f^{lq}(0)\in Y^0,\quad  l=0,1,2\dots.
$$
In this case,  we let $Y^{lq}  = f^{-lq}(Y^0)|\,0$ and consider the unicritical branched covering
$f^q: Y^q\ra Y^0$ of degree $d$. By slight ``thickening'' of the domain of this map (see \cite{M}), 
it can be turned into a unicritical GPL map called the {\it (satellite) renormalization} $Rf$ of $f$.%
\footnote{In the context of GPL maps we use the term ``DH renormalization''
to distinguish it from the generalized renormalization.
In the polynomial case, we refer to it as just ``renormalization'',
as it should not lead to confusion.} 

\ssk
In the satellite renormalizable case, $f$ does not originate any dynasty. 
Otherwise, 
there exists an $l\in \N$ such that $f^{lq}(0)$ belongs
to some puzzle piece $Z_j^1$.
In this case, we let $V^0= f^{-lq}(Z_j^1)|0$ be the {\it first kingdom}, 
and we let $g\equiv g_0: \cup W_i^0\ra V^0$ be the first return map to $V^0$. 
It is easy to check that  $W_i^0 \Subset V^0$. 
%
Let $G_0$ be the associated kingdom map. It originates the dynasty $(g_n, G_n) $ associated with $f$.

\ssk
The map $f$ is called {\it primitively renormalizable} if its dynasty contains  a quadratic-like map $g_n: W^n\ra V^n$
with connected Julia set. This quadratic-like map  is called the {\it (primitive) renormalization} $Rf$ of $f$.
In this case, we cannot construct the next kingdom map  $G_n$,
so the dynasty terminates. 
It also terminates if the map $g$ is combinatorially non-recurrent. 
Otherwise, the process can be continued indefinitely,
and the dynasty $(g_n, G_n)$ is eternal.

If the map $f$ is renormalizable (either in the satellite or in the  primitive sense), 
we can take its renormalization $Rf$
and consider its dynasty. If $Rf$ is renormalizable, 
we can pass to the second renormalization $R^2 f$,   and so on. 
If the map $f$ is at most finitely renormalizable, 
in the end we obtain a non-renormalizable quadratic-like map $R^m f$.
This is the map we will be working with. 
So, {\it in what follows we will assume that $f$ itself is non-renormalizable.}

From now on, we can forget about the original polynomial $f: z\mapsto z^d +c $,
and  replace it with the first map $g: \cup W_i^0\ra V^0$ of the associated dynasty.

\subsection{Extensions}

Let us begin with a trivial but  useful observation:

\begin{lem}[Telescope]\label{telescope}
  Let $X_k$ be a sequence of topological disks, $k=0,1,\dots, m$,
and let $\phi_k : X_k \ra \phi(X^k)$  be  branched coverings of degree $d_k$ such that
$\phi(X_k)\supset X_{k+1}$. Let $\Phi= \phi_{n-1}\circ\dots\circ \phi_0$ (wherever it is defined),
and let $P\subset X_0$ be a component of its domain of definition.
Then $\Phi: P \ra V_n$ is a branched covering of degree at most $d_0\cdot\cdot\cdot d_{n-1}$.
\end{lem}

\begin{lem}\label{landing extension}
  Let $g^m z\in A $ be the first landing  of the $\orb(z)$ at $A$.
Then there exists a puzzle piece $P\ni z$ such that $g^m$ univalently maps $P$ onto $U $. 
\end{lem}  

\begin{pf}
  Let $P= g^{-m} (U) |z$. Then $g^m = h^k \circ g^s$, where $g^s(z)$ is the first landing of $\orb z$
at $U$, $h=h_U: \cup B_i\ra U$ is the generalized renormalization on $U$, and $k$ is the first landing moment
of $\orb_h(g^s(z))$ at $A$. It is a simple exercise to show that $g^s$ is univalent
on $g^{-s}(U)|z$. Moreover,  $h$  univalently maps each non-central component $B_i$, $i>0$, onto $U$.
Now the first assertion follows from the Telescope Lemma.

\end{pf}

\begin{cor}\label{return extension}
  Let $z\in A$, and let $g^m z\in A$ be the first return of the $\orb(z)$ to $A$.
Let  $P = g^{-m}(A)| z$.                     
If  $P$ is not critical then the  map  $g^m: P \ra A$   is univalent.
 Otherwise  $g^m: P \ra A$ is  a unicritical branched covering of  degree $d$. 
\end{cor}

\begin{pf}
 Decompose  $g^m: P\ra A$ as $g: P\ra g(P)$ and the first landing map $g^{m-1}: g(P)\ra A$.
\end{pf}

\comm{***********************************************************
\ssk\nin {\rm (ii)}
  Let  $z=0$, and let $D$ be any component of $g^{-m}(U)$. 
Then  $g^m: D\ra U$ has at most one critical point. 
\end{cor}  

\begin{pf}
(i) The first assertion follows from the fact that at the first moment $T$ when $g^T A=W\ni 0$,
    we have: $g^T A \supset U$. 
    The rest follows from  Lemma \ref{landing extension} by decomposing $g^m|\sm P$ as $L\circ g$,
    where $L$ is the first landing map on $g(P)$. 

\ssk
(ii)
Consider the orbit of $D$,   
$$
        D_k= g^k D, \quad k=0,1,\dots, m-1.
$$
If there are no critical pieces in it, then the map $g^m: D\ra U$ is univalent.
Otherwise, let $D_s$ be the last critical piece in this orbit.  Then $g^{m-s}: D_s\ra U$ is unicritical. 
Moreover, by (i) $D_s\subset A$, and hence the first return time of the critical point  to $D_s$ is at least $m$.
It follows that the preceding  puzzle pieces, $D_0,\dots, D_{s-1}$, are off-critical,
 and the conclusion follows.   
\end{pf}
**********************************************************}

Applying this to the first return of the critical point to $A^{n-1}=V^n$, we obtain:

\begin{cor}\label{space around W}
   The map $g_n: W^n \ra V^n$ admits an analytic extension to a puzzle piece $\tl W^n\Supset W^n$
such that $g_n: \tl W^n \ra U^{n-1}$ is a unicritical branched covering of degree $d$.  
Moreover, $\tl W^n\subset V^n$.
\end{cor}

Let us now construct similar extensions for kingdom maps: 

\begin{lem}\label{cascade extension}
 There is puzzle piece $\tl A \Supset A$ such that  
the map $G: A \ra W $  admits a unicritical degree $d$  extension to a map $\tl A  \ra V$.
Moreover, $\tl A\subset \De$ where $\De$ is the king's apartment.
\end{lem}

\begin{pf}
The map $G: A\ra W$ can be decomposed as $g^k\circ  g^N$ where $N$ is the length of the central cascade of $g$,
 and $k$ is the first entry time of $\orb_g (g^N(A))$ to $W$ (recall that $g^N(A)\subset V\sm W$). 
The map $g^N: A\ra g^N(A)$ admits an analytic  extension to  a unicritical $d$-to-1 branched covering
$g^N: \De\ra W_j$ for some $j>0$, while $g: W_i\ra V$ is a biholomorphic isomorphism for any $i>0$.
Now the conclusion follows by the Telescope Lemma. 
\end{pf}


Let us define {\it enlargements} $\hat E^i$ of domains $E^i$
of the Modified Principal Nest (\ref{E})  as follows: 
$\hat W^n =V^n$ and $\hat V^n = \De^{n-1}$.
We also have the {\it buffers} $\tl E^i\subset \hat E^i$
 constructed in Corollary  \ref{space around W} and Lemma \ref{cascade extension}.
These lemmas tell us that  any map $\psi_i$ analytically extends to a unicritical  $d$-to-1 branched covering 
$\psi_i: \tl E^i\ra \hat E^{i-1}$. 
 For $i<k$, let 
$$
  \Phi_{i,k} = \psi_{i+1}\circ\dots\circ \psi_k: E^k\ra E^i.
$$

By the Telescope Lemma, we have:

\begin{lem}\label{Phi}
    For $0<i<k$, the map $\Phi_{i,k}$ admits an analytic extension to a $d^{k-i}$-to-1 branched covering
from some puzzle piece $F^k\supset E^k$ onto $\hat E^i$.
\end{lem}

\subsection{Travel times} 
  Consider two puzzle pieces $P$ and $Q$ for some GPL or kingdom map $F$.
If $F^l P = Q$, we let $\Time_F (P,Q) =l$
(note that time $l$ is uniquely determined). 
For the ``absolute time'' measured with respect to the initial map $g$, we use notation
 $\Time(P,Q) \equiv \Time_g (P,Q)$.

Let

\ssk \nin $\bullet$ $T_n= \Time (A^n, W^n)$, 
i.e., $ G |A^n = g^{T_n}|A^n$
(the travel time that the king  spends  away from his castle);

\ssk\nin $\bullet$ $t_n= \Time(W^n, V^n)$, i.e., $g_n| W^n = g^{t_n}$;  

\ssk\nin $\bullet$ $s_n = \Time (W^n, W^{n-1})  = t_n+ T_{n-1}$ for $n\geq 1$;  $s_0=t_0=\Time (W^0, V^0)$. 


\begin{lem}\label{t-T-s}
The travel times satisfy the following inequalities:
    $$     t_n\geq T_{n-1}\,; \quad     T_n \geq s_n \, ; \quad    s_n \geq 2 s_{n-1} . $$

\end{lem} 

\begin{pf}
By definition, $g_n(W^n)$ is the first return of $W^n$ to $V^n=A^{n-1}$ under iterates of $G_{n-1}$,
so that $g_n|W^n = G_{n-1}^k|W^n$ for some $k\geq 1$.
Hence 
$$
   g^{t_n}|W^n = g_n|W^n = G_{n-1}^{\circ(k-1)}\circ G_{n-1}|W^n = g^s\circ g^{T_{n-1}}|W^n
$$ 
for some $s\geq 0$,
and the first inequality follows. 

For the second inequality, notice that $T_n$ is the first return time of the critical orbit
to $W^n$ after the first entry to the annulus $V^n\sm W^n$. 
The first entry to $V^n\sm W^n$ occurs at time $\geq t_n$ (since $t_n$ is the first return time of 0 to $V^n$).
Return back to $W^n$ from $V^n\sm W^n$ occurs  at time $\geq T_{n-1}$ (since $T_{n-1}$ is the first moment  $T$
when $f^T(V^n)\cap V^n\not=\emptyset$).

Now the  third inequality follows:
$$
  s_n= t_n+T_{n-1} \geq 2 T_{n-1} \geq  2 s_{n-1}. 
$$

\end{pf} 

\begin{cor}\label{growth}
For any $g$  
we have:   $t_n=\Time(W^n,V^n) \geq 2^{n-1}.$
\end{cor}

\begin{lem}\label{Time}
    $\Time (W^n,W^{n-2}) \geq \Time (V^n,V^0). $
\end{lem}

\begin{pf}
  We have:   
$$
  \Time (W^n,W^{n-2}) = s_n+s_{n-1}= t_n+T_{n-1}+s_{n-1},
$$
while
$$
   \Time (V^n,V^0) = \Time(A^{n-1}, W^{n-1})+ \Time(W^{n-1},V^0) =
$$ 
$$ 
   = T_{n-1} + s_{n-1}+\dots + s_0 .
$$ 
   Thus, the desired inequality is reduced to:
$$
  t_n \geq s_{n-2}+\dots + s_0.
$$
Now the first two inequalities of Lemma \ref{t-T-s} imply that $t_n\geq s_{n-1}$,
and the last one implies that $s_{n-1} \geq s_{n-2}+\dots + s_0 .$

\end{pf}

  Take some $W^n$, and let $l_0$ be the smallest $l\geq \Time(V^n, V^0)$ such that
$g^l (W^n)\subset W^0$. 

\begin{lem}\label{l}
  $l_0\leq \Time(W^n, W^{n-2})$.
\end{lem}

\begin{pf}
  Let $p=\Time(V^n, V^0)$, $l=\Time(W^n, W^{n-2})$. 
By Lemma \ref{Time}, $l \geq p$. Moreover, $g^l (W^n)= W^{n-2}\subset W^0$.
Hence $l\geq l_0$ by definition of $l_0$. 
\end{pf}

We will now make some combinatorial choices.

Fix some (big) $m$. 
Let  $l_0< l_1< l_2<\dots<l_m $ be the $m$ consecutive return moments of the $\orb W^n$ to $W^0$.
In other words, 
$$
   g^{l_k}(W^n) = h^k (g^{l_0}(W^n)),
$$ 
where $h$ is the generalized renormalization of $g$  on $W^0$.  

Let $n> \log_2 m +5$. 
Then by Corollary \ref{growth}, 
\begin{equation}\label{rel}
       \Time_h (W^{n-2}, V^{n-2}) >   \Time_h (W^{n-3}, V^{n-3}) \geq 
\end{equation}
$$
                \geq \Time_{g_1} (W^{n-3}, V^{n-3}) >  2^{n-5} > m.
$$
or, in the absolute time: 

\begin{equation}\label{absolute}
   l_m - l_0 < \Time(W^{n-3}, V^{n-3}) <  \Time (W^{n-2}, V^{n-2}). 
\end{equation}

Putting this estimate together with  Lemma \ref{l}, we conclude:

\begin{lem}\label{last moment}
  $l_m < \Time (W^n, V^{n-2})$.
\end{lem}

\subsection{Degrees}
  
Let  $O= (g^{l_k} (W^n))_{k=0}^m$.
By Lemma \ref{last moment}, $O$ is contained in  the piece $\TT$ of  $\orb_h W^n$
beginning with $W^n$ and ending with  $V^{n-2}$. 
Let us split $\TT$ into five pieces.
 Namely, let $\TT_i$  be the pieces of $\TT$ between
two consecutive domains, $E^i$ and $E^{i-1}$, of the sequence  
\begin{equation}\label{domains}
                W^n \equiv  E^{2n+1}, E^{2n}, \dots , E^{2n-4} \equiv  V^{n-2}.
\end{equation}
Let $O_i= \TT_i\cap O$.

By (\ref{rel}), each $\TT_i$ has length bigger than  $m$.
Hence at most two of the pieces $O_i$ are non-empty, and so one of them
contains at least $m/2$ elements. Let now $O_i$ stand for such a piece.


Let us consider the enlargement $\hat E^{i-1}$ of $E^{i-1}$.
Notice that it is  contained in $W^{n-3}$. 
Let us pull  $\hat E^{i-1}$ back along the $h$-orbit of $W^n$.
It inscribes every domain of this orbit,  $W^n, h(W^n), \dots, h^s(W^n)= E^{i-1}$, 
into  a bigger {\it buffer} domain $F, h(F),\dots, h^s(F) = \hat E^{i-1}$. 

By Lemma \ref{Phi}, we have:

\begin{lem}\label{degree bound}
   The map $h^s: F \ra \hat E^{i-1}$ has degree at most $d^5$.
\end{lem}

Moreover, 

\begin{lem}\label{disjoint}
   The domains $h^k(F)$ enclosing the domains of $O_i$ are pairwise disjoint.
\end{lem}

\begin{pf}
  Otherwise there would be two nested domains $h^k(F) \subset h^s (F)$, $k<s$. 
Let $L = s-k$. Pushing $h^k(F)$ forward to $\hat E^{i-1}$ we see that
$h^L (\hat E^{i-1})\supset \hat E^{i-1}$. All the more, $h^L (W^{n-3}) \supset W^{n-3}$,
so that $L \geq  \Time_h (W^{n-3}, V^{n-3})$.

On the other hand, by  (\ref{absolute}),
 $$
    L < \Time_h (W^{n-3}, V^{n-3}), 
$$
contradiction.
\end{pf}

Let us now consider some domain $\La= g^{l_k}(W^n)\in O_i$, and let $\La' = g^{l_k} (F)$ be its  buffer.
Since there is a biholomorphic push-forward $(\La', \La)\ra ( \hat E^{i-1}, E^{i-1} )$, we have:

\begin{lem}\label{moduli of buffers}
   $ \mod(\La'\sm \La) = \mod (\hat E^{i-1}\sm  E^{i-1}).$
\end{lem} 
  
Let $\Upsilon= g^{-l_k}(V^0)|0$.

\begin{lem}\label{outer degree}
  We have: $W^n\subset \Upsilon \subset V^n$ and 
$$
   \deg(g^{l_k}: \Upsilon\to V^0) \leq d^{2n+m}.
$$
\end{lem}

\begin{pf}
  The first inclusion is trivial. The  second inclusion, $\Upsilon\subset V^n$,
follows from  $l_k\geq l_0\geq \Time(V^n, V^0)$. 

Let us estimate the degree. Let $s= \Time(V^n,V^0)$. 
Then  
$$
     \deg(g^s : V^n \ra V^0)=d^{2n}. 
$$

Let us now consider the first landing map $H$ to $W^0$.
It is easy to see that each component $Q_j$ of the domain of $H$ 
is mapped biholomorphically onto  $W^0$ and, moreover,
 $H:  Q_j\ra W^0$  admits an extension
to a biholomorphic isomorphism $\tl Q_j\ra V^0$. 
Let $\Upsilon_i = g^{l_i}(\Upsilon) $.  Then we have:
$$
  \Upsilon_0  =       H \circ g^s (\Upsilon)
$$
and 
$$
     \Upsilon_{i+1} = H\circ (g | W^0)|\,  \Upsilon_i , \quad i=0,1,\dots, k-1 \leq m-1, 
$$
and the Telescope Lemma concludes the proof. 
\end{pf}

\subsection{Summary}\label{summary} 
We fix an arbitrary $m$ and let $n>\log_2 m +5$. 
Then for any domain $\La =\La_k = g^{l_k}(W^n) \in O_i$, 
the map $\Psi=\Psi_k= g^{l_k}: W^n \ra \La$ admits a holomorphic extension to a branched covering
\begin{equation}\label{3 domains}
    \Psi: (\Upsilon, F, W^n) \ra (V^0, \La', \La)
\end{equation}
such that:

\begin{itemize}

\item [(P1)] $\deg(\Psi: \Upsilon\ra V^0)\leq d^{2n+m}$;

\item[(P2)] $\deg(\Psi: F\ra \La')\leq d^5$;

\end{itemize}

\nin and

\begin{itemize}
\item [(P3)] $\Upsilon\subset V^n$;   

\item[(P4)] $\mod(\La' \sm \La) = \mod (\hat E^{i-1}\sm  E^{i-1})$.
\end{itemize}

\ssk\nin
Moreover, there are at least $m/2$ domains $\La_k$ in the orbit $O_i$, and their buffers $\La_k'$ are 
pairwise disjoint.

\comm{

\begin{lem}\label{once}
All the domains $g^{l_k} (W^n)$, $k=l_0, \dots, m$,
except at most one, are disjoint from  $W^{n-3}$.
\end{lem}

\begin{pf}
  By Lemma \ref{last moment},  all these  domains have depth bigger than the depth of $V^{n-2}$,
all the more, bigger than the depth of $W^{n-3}$.
Hence if some $g^{l_k}(W^n)$ intersects $W^{n-3}$ then it is contained in $W^{n-3}$.
But then by (\ref{absolute}), all the further domains are contained in 
$$ 
  \bigcup_{s=1}^{t_{n-3}-1} g^s (W^{n-3}),
$$
which is disjoint from $W^{n-3}$. 
\end{pf}
end comm}

\section{Quasi-Additivity Law and Covering Lemma}\label{covering lemma sec}

\proclaim Quasi-Additivity Law (\cite{covering lemma}, {\rm \S 2.9}). 
    Fix some $\eta>0$. 
Let $W\Subset V$ and  $\La_i\Subset \La_i'\Subset W $, $i=1,\dots, m$, 
be topological disks such that the closures of $\La_i'$ are pairwise disjoint. 
Then there exists a $\de_0 > 0$ (depending on $\eta$ and $m$) such that: \\
If  for some $\de\in (0, \de_0)$,
$\mod (V\sm \La_i) < \de$ while $\mod(\La_i'\sm \La_i) > \eta\de$, then
$$
  \mod(V\sm W) < \frac{C \eta^{-1} \de}{m},
$$
where $C$ is an absolute constant.

\proclaim Quasi-Invariance Law/Covering Lemma \cite{covering lemma}. 
Fix some  $\eta>0$. 
 Let $U\Supset A'\Supset A$ and $V\Supset B' \Supset B$ be two nests of topological disks.
Let $g: (U, A', A) \ra (V, B', B)$ be a branched covering between the respective disks.
Let $d = \deg(A'\ra B')$ and  $D = \deg(U\ra V)$.
Assume  
$$
   \mod(B'\sm B) > \eta \mod(U\sm A).
$$
If $\mod(U\sm A) < \eps(\eta, D) $ then 
$$
      \mod (V\sm B)  <  C\eta^{-1} d^2 \mod(U\sm A),
$$
where $C$ is an absolute constant.

\section{A priori bounds}

The following  lemma tells us that if some principal  modulus is very small then it
should be even smaller on some preceding level:

\begin{lem}\label{growth-lem}
There exist $n=n(d)\in \N$ and $\eps=\eps(d,n)>0$  such that: 
If on some level $q\geq n$, $  \mod( V^q\sm W^q ) <\eps $, 
then on some previous level $p<q$ we have:
\begin{equation}\label{decrease}
    \mod  (V^p \sm W^p) <  \frac{1}{2} \mod( V^q\sm W^q ).   
\end{equation}
\end{lem}

\begin{pf}
We will use the set-up of \S \ref{summary}, except that the base GPL map $g$ will not be $g_0$
but rather $g_s$ on some deeper level.
Let us fix some $m > 16 C^3 d^{23}$,
where $C$ is the maximum of  constants in the First and Second Covering Lemmas.
Let $q> n > \log_2 m+ 5$. 
We take $g = g_{q-n}$ as the base map
and consider the associated 3-domain branched covering $\Psi=\Psi_k$ (\ref{3 domains}) 
$$
   \Psi: (\Upsilon, F, W^q) \ra (V^{q-n}, \La', \La),  
$$
where $\La=\La_k$ is one of the domains of the orbit $O_i$.  
Set $\eta=1/2d$ for the First Covering Lemma. 
Let us consider two cases:

\msk {\it Case 1.} Assume that for some domain $\La\in O_i$, 
$$
  \mod (\La'\sm \La)< \frac{1}{2d} \mod (V^q\sm W^q).
$$
         By Property (P4), 
$\mod(\La' \sm \La) = \mod (\hat E^{i-1}\sm  E^{i-1})$, which is equal to either
$\mod(V^{(i-2)/2}\sm W^{(i-2)/2})$ (if $i$ is even) or to
$$
   \mod(\De^{(i-3)/2}\sm A^{(i-3)/2}) \geq \frac{1}{d} \mod (V^{(i-3)/2}\sm W^{(i-3)/2})  \quad \mbox{if $i$ is odd. }
$$   
In both cases we conclude that (\ref{decrease}) holds  for   
$p$ which  is equal  to either $(i-2)/2$ or $(i-3)/2$.
(Note that $p<q$ since by construction of the  buffers, $i-1 <  2q+1$.)     

\msk{\it Case 2.} Assume that for all $\La_k\in O_i$, 
\begin{equation}\label{1/2}
  \mod (\La_k'\sm \La_k) \geq \frac{1}{2d} \mod (V^q\sm W^q).
\end{equation}
Then  the  Covering Lemma is applicable to every map $\Psi = \Psi_k$, 
provided $\eps=\eps(d,n)$ is sufficiently small. 
It yields:
\begin{equation}\label{outer mod}
  \mod (V^{q-n}\sm \La_k) \leq 2C d^{11} \mod(\Upsilon\sm W^q) \leq 2C d^{11} \mod(V^q\sm W^q). 
\end{equation}
Estimates (\ref{1/2}) and (\ref{outer mod}) show that the Quasi-Additivity Law is applicable
with $\eta=1/4Cd^{12}$.
 Since there are at least $m$ domains $\La_k\subset \La_k'\subset W^{n-q}$ in the orbit $O_i$,
it implies:
$$
      \mod(V^{n-q}\sm W^{n-q}) \leq \frac{8 C^3 d^{23}}{m} \mod(V^q\sm W^q) < \frac{1}{2} \mod(V^q\sm W^q),
$$
and we are done.  
\end{pf}

Lemma \ref{growth-lem} immediately yields Theorem B from the Introduction.

\end{document}